
  \font\scX=cmcsc10
  \font\rmVIII=cmr8
  \font\rmVI=cmr6
  
  \font\bfXII=cmbx10 scaled \magstep1
  \def\Stratila{Str\v atil\v a}
  \def\:{\colon}
  \def\cstar{$C^*$}
  \def\arw{\rightarrow}
  \def\H{{\cal H}}
  \def\Z{{\bf Z}}
  \def\N{{\bf N}}
  \def\S1{S^1}
  \def\for{\hbox{\ \ \ for\ }}
  \def\th{\null^{\hbox{\rmVI th}}}
  \def\max{\hbox{\rmVI max}}
  \def\e{{\rm e}}
  \def\sameauthor{\underbar{\hbox{\ \ \ \ \ \ \ \ \ \ }}}
  \def\statement#1#2{\medbreak\noindent{\bf#1.\enspace}{\sl#2\par}
\ifdim\lastskip<\medskipamount\removelastskip\penalty55 \medskip\fi}
  \def\proof{\smallskip \noindent {\bf Proof. }}
  \def\square{\hbox{$\sqcap\!\!\!\!\sqcup$}}
  \def\endproof{\hfill \square}
  \def\cite#1{{\rm[\bf #1\rm]}}
  \def\stdbib#1#2#3#4#5#6{\smallskip
  \item{[#1]} #2, ``#3,'' {\sl #4} {\bf #5} #6.}
  \def\bib#1#2#3#4{\smallskip
  \item{[#1]} #2, ``#3,'' {#4.}}

  \def\Bratteli{1}
  \def\EL{2}
  \def\NP{3}
  \def\Loring{4}
  \def\PV{5}
  \def\PVEmbedd{6}
  \def\Ped{7}
  \def\Power{8}
  \def\Putnam{9}
  \def\SV{10}

  \magnification=\magstep1
  \nopagenumbers
  \voffset=2\baselineskip
  \advance\vsize by -\voffset
  \headline{\ifnum\pageno=1
    \hfil
  \else
    \ifodd\pageno
      \scX\hfil af-algebras and partial automorphisms\hfil\folio
    \else
      \scX\folio\hfil ruy exel\hfil
    \fi
  \fi}

  \null
  \vskip 1truein

  \baselineskip=18pt
  \centerline{\bfXII APPROXIMATELY FINITE C\raise.5ex\hbox{*}-ALGEBRAS
AND}
  \centerline{\bfXII PARTIAL AUTOMORPHISMS}
  \centerline{(Preliminary version)}

  \bigskip
  \bigskip
  \centerline{\scX Ruy Exel\footnote
  {\raise.5ex\hbox{*}}
  {\rmVIII On leave from the University of S\~ao Paulo.}}

  \bigskip
  \baselineskip=12pt
  \centerline{\it Department of Mathematics and Statistics}
  \centerline{\it University of New Mexico}
  \centerline{\it Albuquerque, New Mexico 87131}
  \centerline{e-mail: exel@math.unm.edu}

  \vfill
  \midinsert\narrower\narrower\noindent
  {\bf Abstract.} We prove that every AF-algebra is isomorphic to a
crossed product of a commutative AF-algebra by a partial automorphism.
The case of UHF-algebras is treated in detail.
  \endinsert\vfill\eject

  \beginsection 1. Introduction

Let $\alpha$ be an action of the group of integers on a unital
\cstar-algebra $A$.  It is an easy consequence of the
Pimsner--Voiculescu exact sequence \cite{\PV}, that the $K_1$-group of
the crossed product algebra $A\times_\alpha \Z$ is \underbar{never}
trivial (at least if the unit of $A$ is non-zero in $K_0(A)$).  On the
other hand, one of the important properties of approximately finitely
\cstar-algebras (AF-algebras for short) is that their $K_1$-groups are
\underbar{always} trivial.  This says that AF-algebras are in a class
apart from algebras obtained by crossed products by \Z, at least from
the point of view of $K$-theory.

  Despite this fact, many authors have been extremely successful in
finding deep relationships between these seemingly different classes
of \cstar-algebras.  The first such interconnection that comes to mind
is, of course the famous embedding of the irrational rotation
\cstar-algebras \cite{\PVEmbedd} into AF-algebras, obtained by Pimsner
and Voiculescu.  The work of Loring \cite{\Loring} as well as that of
Elliott and Loring \cite{\EL} on AF-embeddings of the algebra of
continuous functions on the two-torus, extend that relationship even
further.

  \Stratila\ and Voiculescu have characterized AF-algebras as concrete
representations of certain \cstar-dynamical systems in (I.1.10) of
\cite{\SV}, thus adding to the evidence that AF-algebras are a sort of
a crossed product.  The work of Putnam \cite{\Putnam}, is also rather
intriguing as he shows that certain crossed product \cstar-algebras
associated to minimal homeomorphisms of the Cantor set contain, as
well as are contained in, AF-algebras.

  The theory of inductive limits of ``circle algebras'' developed by
Elliot shows one of the strongest evidences to date, of the close
connections between crossed products and inductive limits: the
irrational rotation \cstar-algebras, themselves, can be described as
inductive limits of circle algebras.

  More recently, we have seen substantial progress in the direction of
extending both the concept of AF-algebras, as well as that of crossed
products, to the realm of non-self-adjoint operator algebras. A good
source for references for this is Power's notes on subalgebras of
\cstar-algebras \cite{\Power}.

  The cases mentioned above certainly do not exhaust the list of
interesting connections between crossed products and AF-algebras, but
they do come close to exhausting the authors knowledge of the
literature on the subject of AF-algebras, were we claim no expertise.

  Nevertheless, the main result of the present work, Theorem (2.6),
describes what we believe is the most straightforward possible
relationship between AF-algebras and crossed products.  Namely, that
all AF-algebras are isomorphic to a crossed product of an AF-masa
(maximal abelian subalgebra), not by an automorphism, but by a partial
automorphism. We feel that this fact may well be one of the bases for
the the rich parallel between the two classes of \cstar-algebras under
consideration.

  Partial automorphisms were introduced by the author in \cite{\NP}
were the reader will find the main definitions and results on
generalized crossed products, as well as the facts on actions of the
circle that we shall need here.

  This work is organized in two sections that follow the present
introduction. In section (2) we study circle actions and covariant
homomorphisms on finite dimensional \cstar-algebras, since these are
the building blocks for AF-algebras. Inductive limits are then
considered and actions on the building blocks are extended to the
limit algebra. A somewhat canonical circle action is then defined on
every AF-algebra and our main result is derived from Theorem (4.21) in
\cite{\NP}.

  In section (3) we consider the example of uniformly hyperfinite
\cstar-algebras (UHF-algebras for short) and we obtain a result very
close to Putnam's \cite{\Putnam} study of crossed products of the
Cantor set by ``odometer'' maps.  Since we are not restricted to
classical crossed products, we may describe any UHF-algebra, not only
as a subalgebra, but as the whole crossed product \cstar-algebra of
the Cantor set by a partially defined ``odometer'' map which also
plays a role in \cite{\Putnam}.

  Our treatment of the example of UHF-algebras is only slightly
original since it is basically contained in Putnam's work.
Nevertheless, the unifying view of AF-algebras provided by partial
automorphisms, sheds a new light on these examples.

  The question of which crossed products by automorphisms yield
AF-algebras, which so far had an obvious negative answer, suddenly
gains relevance when one thinks of partial, instead of globally
defined automorphisms.  If properly interpreted, an answer to this
question is provided in section (3) of \cite{\Putnam}. Following the
notation there, the subalgebras $A_Y$ are precisely the crossed
product of $C(X)$ by the partial automorphism obtained by restriction
of the minimal homeomorphism $\varphi$ to the ideal $C_0(X-Y)$.

  As far as notation is concerned, we should warn the reader of a
slightly non-standard convention from \cite{\NP} used here: if $X$ and
$Y$ are subsets of a \cstar-algebra, then $XY$ denotes the
\underbar{closed} linear span of the set of products $xy$ with $x \in
X$ and $y \in Y$.

  \beginsection 2. The Main Result

Let $A$ and $B$ be \cstar-algebras and let $\alpha$ and $\beta$ be
actions of $\S1$ on $A$ and $B$, respectively. If $\phi\: A\arw B$ is
a covariant homomorphism (always assumed to be star preserving), then
it is clear that $\phi(A_n)\subseteq B_n$, where $A_n$ and $B_n$ are
the corresponding spectral subspaces. This obviously implies that
$\phi(A_n^*A_n)\subseteq B_n^*B_n$ for each $n$.

  Let us now suppose that both $\alpha$ and $\beta$ are regular
actions in the sense of (4.4) in \cite{\NP}. That is, there are maps
  \medskip
  \itemitem{(i)} $\theta_A\: A_1^* A_1 \arw A_1 A_1^*$
  \medskip
  \itemitem{(ii)} $\lambda_A\: A^*_1\arw A_1 A_1^*$
  \medskip
  \itemitem{(iii)} $\theta_B\: B_1^* B_1 \arw B_1 B_1^*$
  \medskip
  \itemitem{(iv)} $\lambda_B\: B^*_1\arw B_1 B_1^*$
  \medskip
  \noindent satisfying the conditions described in \cite{\NP}.

  \statement{2.1. Definition}{A covariant homomorphism $\phi\: A\arw
B$ is said to be regular (with respect to a given choice of
$\theta_A$, $\lambda_A$, $\theta_B$ and $\lambda_B$) if for any $x^*$
in $A_1^*$ and $a$ in $A_1^*A_1$ one has
  \medskip
  \itemitem{(i)} $\phi(\lambda_A(x^*))=\lambda_B(\phi(x^*))$
  \medskip
  \itemitem{(ii)} $\phi(\theta_A(a))=\theta_B(\phi(a)).$}

  \medskip Let us now analyze, in detail, an example of a regular
homomorphism in finite dimensions, which will prove to be crucial in
our study of AF-algebras.

  Denote by $M_k$, the algebra of $k\times k$ complex matrices.
  Let $A$ be a finite dimensional \cstar-algebra, so that $A$ is
isomorphic to a direct sum $A=\bigoplus_{i=1}^n M_{p_i}$. In order to
define a circle action on $A$, let for each $k$, $\alpha^k$ denote the
action of $\S1$ on $M_k$ given by
  $$\alpha_z^k(a) = \pmatrix{
  1 \cr
  & z \cr
  && \ddots \cr
  &&& z^{k-1}
  }\ a\ \pmatrix{
  1 \cr
  & z \cr
  && \ddots \cr
  &&& z^{k-1}
  }^{-1}$$
  for any $a$ in $M_k$ and $z$ in $\S1$. Define an action $\alpha_A$
of $\S1$ on $A$ by $\alpha_A = \bigoplus_{i=1}^n \alpha^{p_i}$. This
action will be called the standard action of $\S1$ on $A$. It should
be noted that the standard action depends on the choice of matrix
units for A.

  One can easily check that $\alpha_A$ is semi-saturated (see 4.1 in
\cite{\NP}). Let us show that it is also regular. For each $k$ let
$s_k$ be the partial isometry in $M_k$ given by
  $$s_k=\pmatrix{
  0 \cr
  1 & 0 \cr
    & 1 \cr
    & & \ddots \cr
    & & & 1 & 0
    }$$
  and put $s_A = \bigoplus_{i=1}^n s_{p_i}$. The maps
  $$\lambda_A\: x^* \in A_1^* \mapsto s_A x^* \in A_1 A_1^*$$
  $$\theta_A\: a \in A_1^*A_1 \mapsto s_A a s_A^* \in A_1A_1^*$$
  satisfy the conditions of (4.4) in \cite{\NP} as the reader may
easily verify, which therefore implies the regularity of $\alpha_A$.

  Now let $B=\bigoplus_{i=1}^m M_{q_i}$ be another finite dimensional
\cstar-algebra and let $\alpha_B$, $s_B$, $\lambda_B$ and $\theta_B$
be defined as above.

  A standard homomorphism is any homomorphism $\phi\: A \arw B$ of the
form
  $$\phi(a_1,\ldots,a_n) =
  (\phi_1(a_1,\ldots,a_n), \ldots,
  \phi_m(a_1,\ldots,a_n))$$
  where, for each $k=1,2, \ldots, m$
  $$\phi_k(a_1,a_2,\ldots,a_n) = \pmatrix{
  a_{i_1} \cr
  & a_{i_2} \cr
  && \ddots \cr
  &&& a_{i_r} \cr
  &&&& 0 \cr
  &&&&& \ddots \cr
  &&&&&& 0}$$
  where $r$ may vary with $k$, the indices $i_1,i_2\ldots i_r$ are
chosen from the set $\{1, 2, \ldots, n\}$, while each $a_i$ is in
$M_{p_i}$. Standard homomorphisms are well known to play an important
role in the theory of AF-algebras \cite{\Bratteli}.

  \statement{2.2. Proposition}{A standard homomorphism is covariant
and regular with respect to the respective standard actions.}

  \proof Let $\phi\: A \arw B$ be a standard homomorphism where $A$
and $B$ are direct sum of matrix algebras, as above. That $\phi$ is
covariant follows from the easy fact that spectral subspaces are
mapped accordingly. In order to verify regularity, let $x \in A_1$.
Then
  $$\phi(\lambda_A(x^*)) = \phi(s_Ax^*) = \phi(s_A)\phi(x^*).$$
  Note that $s_B \phi(s_A)^* \phi(s_A) = \phi(s_A)$, since $\phi(s_A)$
can be thought of as a ``restriction'' of $s_B$ to a subspace of the
initial space of $s_B$. Thus
  $$\phi(\lambda_A(x^*)) = s_B\phi(s_A)^*\phi(s_A)\phi(x^*) = s_B
\phi(s_A^*s_Ax^*) = s_B\phi(x^*) = \lambda_B(\phi(x^*)).$$
  On the other hand, for $a$ in $A_1^*A_1$, of the form $a=x^*y$, with
$x,y\in A_1$, we have, by the same reasons given above
  $$\phi(\theta_A(a)) = \phi(s_Ax^*ys_A^*) = s_B\phi(x^*y)s_B^* =
\theta_B(\phi(a)).$$
  This completes the proof. \endproof

  \statement{2.3. Theorem}{Let $(A^k)_{k\in\N}$ be a sequence of
\cstar-algebras and $\phi^k\: A^k \arw A^{k+1}$ be injective
homomorphisms. Suppose each $A^k$ carry a semi-saturated regular
action $\alpha^k$ of $\S1$. Suppose, in addition, that each $\phi^k$
is covariant and regular with respect to a fixed choice of maps
$\lambda^k$ and $\theta^k$, at each level, as in (2.1). Then the
inductive limit \cstar-algebra $A = \lim_k A^k$ admits a
semi-saturated regular action $\alpha$ of $\S1$. Moreover the
canonical inclusions $\psi^k\: A^k \arw A$ are covariant and regular.}

  \proof Let $A'$ be the union of all $\psi^k(A^k)$ in $A$.  For each
$k$ and each $a\in A^k$ define
  $$\alpha_z(\psi^k(a)) = \psi^k(\alpha_z^k(a)) \for z\in \S1.$$
  Since the $\phi^k$ are covariant, it follows that $\alpha_z$ is well
defined as an automorphism of $A'$ and that it extends to an action of
$\S1$ on $A$, with respect to which, each $\psi^k$ is covariant. So it
is clear that $\psi^k(A_n^k) \subseteq A_n$, where the subscript $n$
indicates spectral subspace. The fact that $\alpha$ is semi-saturated,
thus follows immediately.

  We now claim that each $A_n$ is indeed the closure of
$\bigcup_{k=1}^{\infty}\psi^k(A_n^k)$. In fact, let $P_n^k$ and $P_n$
be the $n\th$ spectral projections corresponding to $A^k$ and $A$,
respectively.  Given $a$ in $A_n$, choose $a^k$ in $A^k$ such that
$a=\lim_k \psi^k(a^k)$. Note that $a = P_n(a) =
\lim_k\psi^k(P_n^k(a^k))$ which proves our claim. We then clearly have
that, for any integer $n$, $A_n^*A_n$ is the closure of
$\bigcup_{k=1}^{\infty}\psi^k({A_n^k}^*A_n^k)$.  For each $k\in \N$,
$x \in A_1^k$ and $a \in {A_1^k}^*A_1^k$ define
  $$\lambda(\psi^k(x^*)) = \psi^k(\lambda^k(x^*))$$
  $$\theta(\psi^k(a)) = \psi^k(\theta^k(a)).$$
  Invoking the regularity of $\phi^k$, we see that both $\lambda$ and
$\theta$ are well defined, the domain of $\lambda$ being
  $\bigcup_{k=1}^{\infty}\psi^k({A_1^k}^*)$
  and that of $\theta$ being
  $\bigcup_{k=1}^{\infty}\psi^k({A_1^k}^*A_1^k)$.
  The hypothesis clearly imply that $\psi^k$, $\theta^k$ and
$\lambda^k$ are isometries, so both $\lambda$ and $\theta$ extend to
the closure of their current domains, that is, $A_1^*$ in the case of
$\lambda$ and $A_1^*A_1$ in case of $\theta$.  One now needs to show
that for any $x, y \in A_1$, $a \in A_1^* A_1$ and $b \in A_1 A_1^*$
  \medskip
  \itemitem{(i)} $\lambda (x^* b) = \lambda(x^*) b$
  \medskip
  \itemitem{(ii)} $\lambda (ax^*) = \theta (a) \lambda (x^*)$
  \medskip
  \itemitem{(iii)} $\lambda (x^*)^* \lambda (y^*) = xy^*$
  \medskip
  \itemitem{(iv)} $\lambda (x^*) \lambda (y^{*})^{*} = \theta (x^* y)$
  \medskip \noindent which are the axioms for regular actions (see 4.4
in \cite{\NP}). These identities hold on dense sets, by regularity of
the $\phi^k$, and hence everywhere by continuity. A few remaining
details are left to the reader. \endproof

  \smallskip A simple, yet crucial consequence of (2.2) and (2.3) is
the existence of certain somewhat canonical circular symmetries on
AF-algebras. More precisely we have:

  \statement{2.4. Theorem} {Every AF-algebra possesses a
semi-saturated regular action of the circle group, such that the fixed
point subalgebra is an AF-masa.}

  \proof For the existence of the action it is enough to note that any
AF-algebra can be written as a direct limit of direct sums of matrix
algebras, in which the connecting maps are standard homomorphisms in
the above sense \cite{\Bratteli}.  The fixed point subalgebra is the
inductive limit of the direct sum of the diagonal subalgebras at each
stage, hence it is commutative and AF.  The maximality follows from
(I.1.3) in \cite{\SV}.  \endproof

  \smallskip Of course the action described in (2.4) is not unique as
it depends on the choice of a particular chain of finite dimensional
subalgebras, as well as on the choice of matrix units at each stage.
Nevertheless, the isomorphism class of the fixed point algebra will
always be the same \cite{\Power}. With a little reluctance we, thus
introduce the following:

\statement{2.5. Definition}{Given an AF-algebra $A$, the action of
$\S1$ on $A$ described in the proof of (2.4) will be called the
standard action.}

  \smallskip Combining the existence of standard actions with Theorem
(4.21) in \cite{\NP}, brings us to our main result.

  \statement{2.6. Theorem}{Let $A$ be an AF-algebra. Then there exists
a totally disconnected, locally compact topological space $X$, a pair
of open subsets $U, V \subseteq X$ and a homeomorphism $\varphi \: V
\arw U$ such that the covariance algebra for the induced partial
automorphism $\theta \: C_0(U) \arw C_0(V)$ is isomorphic to $A$.}

  \proof Apply Theorem (4.21) in \cite{\NP} to the action of $\S1$ on
$A$ provided by (2.4). \endproof

  \beginsection 3. An Example

We would now like to analyze, in detail, the case of UHF-algebras.
Given an UHF-algebra $A$, we propose to exhibit, in a very concrete
fashion, the space $X$, the open sets $U$ and $V$ as well as the map
$\varphi$ referred to in (2.6).

  In the most general case, this is not an easy task. The knowledge
that a \cstar-algebra is a covariance algebra for some partial
automorphism, often comes without any hint about what that partial
automorphism should be. Obtaining a concrete description of it can,
sometimes, be a little tricky.  Nevertheless, when the fixed point
subalgebra is commutative, life is a bit easier.

  \statement{3.1. Lemma}{Let $\alpha$ be a semi-saturated regular
action of $\S1$ on a \cstar-algebra $A$ and let $\theta\: A_1^* A_1
\arw A_1A_1^*$ be as in (4.4) of \cite{\NP}.  Suppose that the fixed
point subalgebra $A_0$ is commutative. Then for every $a\in A_1^* A_1$
of the form $a=x^*y$, where $x,y \in A_1$, one has $\theta(a) = yx^*$.
Therefore $\theta$ is uniquely determined.}

  \proof If we fix a faithful representation of $A$, we may assume
that $A$ is a subalgebra of $B(\H)$, for some Hilbert space $\H$. From
(5.6) in \cite{\NP} we conclude that there is a partial isometry $u
\in B(\H)$ such that $\theta(a) = u a u^*$ for any $a$ in $A_1^*A_1$.
Moreover we have $A_1 = A_1A_1^*u$. Therefore $x=hu$ and $y=ku$ for
some $h, k\in A_1A_1^*$. We have
  $$\theta(a) = uau^* = uu^*h^*kuu^* = h^*k$$
  while
  $$yx^* = k u u^* h^* = k h^*.$$
  The conclusion follows from the commutativity of $A_0$. \endproof

  \smallskip In order to establish our notation let's now carefully
define UHF-algebras (see 6.4 in \cite{\Ped}). Let
$(n_i)_{i=0}^{\infty}$ be a sequence of positive integers. For $k\geq
0$ let $p_k=\prod_{i=0}^{k-1} n_i$ with the convention that $p_0=1$.

  If $\beta=(\beta_0,\ldots,\beta_{k-1})$ is a finite sequence of
integers with $0 \leq \beta_i < n_i$, let $\e_\beta$ be the $p_k\times
p_k$ complex matrix
  $$\e_\beta = \pmatrix{
  0 \cr
  & \ddots \cr
  && 1 \cr
  &&& \ddots \cr
  &&&& 0}$$
  where the only non-zero entry appears in the position
  $$j(\beta) = \beta_0 p_0 + \cdots + \beta_{k-1} p_{k-1},$$
  counting the upper left hand corner as position zero.

  Alternatively, the location of the non-zero entry in $\e_\beta$ can
be determined as follows: divide the above matrix in $n_{k-1}\times
n_{k-1}$ blocks of size $p_{k-1}\times p_{k-1}$ each. Then, counting
down along the diagonal, from zero to $\beta_{k-1}$, locate a block.
Repeat this process within that block until a $1\times 1$ block is
identified. There you will find the entry ``1''.

  \smallskip
  Consider the embedding $\phi_k\: M_{p_k} \arw M_{p_{k+1}}$ given by
  $$\phi_k(a) = \pmatrix{
  a \cr
  & a \cr
  && \ddots \cr
  &&& a}$$
  where each block ``$a$'' appears $n_k$ times.

  The UHF-algebra associated to the supernatural number
$\prod_{i=0}^{\infty} n_i$ is the inductive limit of the sequence
  $$M_{p_0} {\buildrel \phi_0 \over \longrightarrow}
  M_{p_1} {\buildrel \phi_1 \over \longrightarrow}
  M_{p_2} {\buildrel \phi_2 \over \longrightarrow} \cdots$$
  Denote the resulting \cstar-algebra by $A$. The standard action of
$\S1$ on $A$ has, as its fixed point algebra, the subalgebra $A_0$ of
$A$, generated by the $\e_\beta$ for all finite sequences $\beta$ as
above.  The spectrum of $A_0$ is the compact space
  $$\sigma(A_0) = \prod_{i=0}^{\infty} \{0, 1, \ldots, n_i-1\}$$
  which is well known to be homeomorphic to the Cantor set.  Each
$\e_\beta$ corresponds to the characteristic function of the clopen
set formed by all infinite sequences in $\sigma(A_0)$ containing
$\beta$ as the initial segment.

  We shall now determine the ideals $A_1^*A_1$ and $A_1A_1^*$ of the
fixed point algebra $A_0$ as well as the partial automorphism
$\theta\: A_1^*A_1 \arw A_1A_1^*$ referred to in (2.6). Reasoning as
in the proof of (2.3), one sees that the first spectral subspace
$A_1$, is linearly generated (in the sense of closed linear span) by
the matrix units $\e_{j+1,j}$ in $M_{p_k}$ for all $k$ and for $0\leq
j \leq p_k-2$. Therefore it is clear that $A_1^*A_1$ is linearly
generated by the projections $\e_{jj}$ in $M_{p_k}$ for $0\leq j \leq
p_k-2$. This corresponds exactly to the projections $\e_\beta$ where
$\beta$ is not equal to
  $$\beta_{\max}^k := (n_0-1,n_1-1, \ldots, n_{k-1}-1).$$
  A similar analysis shows that $A_1A_1^*$ is generated by the
$\e_\beta$ with $\beta$ not equal to
  $${\bf 0}^k :=(0, 0, \ldots, 0).$$
  It follows that $A_1^*A_1 = C_0(U)$ and $A_1A_1^* = C_0(V)$ where
$U$ and $V$ are the open subsets of $\sigma(A_0)$ given by
  $U=\sigma(A_0)-\{\beta_{\max}\}$ and
  $V=\sigma(A_0)-\{{\bf 0}\}$
  where $\beta_{\max}$ and ${\bf 0}$ are the infinite sequences
  $(n_0-1,n_1-1, \ldots)$ and $(0, 0, \ldots)$, respectively.

  In order to find an expression for the partial automorphism
$\theta\: C_0(U) \arw C_0(V)$, we use (3.1). Let $\beta =
(\beta_0,\ldots,\beta_{k-1}) \neq \beta_{\max}^k$ so that $\e_\beta$
is in $C_0(U)$. Also note that $j(\beta) = \beta_0 p_0 + \cdots +
\beta_{k-1} p_{k-1}$ is strictly less than $p_k-1$. Let
$x=\e_{j(\beta)+1,j(\beta)} \in M_{p_k}$ so that $x$ is in $A_1$. Then
$\e_\beta = x^*x$ so, by (3.1) we have
  $$\theta(\e_\beta) = x x^* = \e_{j(\beta)+1,j(\beta)+1} =
\e_{\beta'},$$
  where $\beta'$ is such that $j(\beta') = j(\beta) + 1$. To obtain
$\beta'$ one should interpret $\beta$ as a sequence of digits and add
$(1,0, \ldots, 0)$ to $\beta$, with carry over to the right.

  The fact that the $\e_\beta$ generate $A_0$ implies that the map
$\varphi\:V \arw U$ inducing $\theta$, must be the inverse of the map
{}from $U$ to $V$, obtained by addition of the infinite sequence $(1,
0, \ldots)$, again with carry over to the right. This map is
sometimes referred to as the ``odometer'' map for obvious reasons.

  Note that, since $\beta_{\max}$ is not in the domain of definition
of $\varphi^{-1}$, the process of carrying over will always terminate
in finite time. We have then proved the following:

\statement{3.2. Theorem}{Let $(n_i)_{i=0}^{\infty}$ be a sequence of
positive integers and let
  $$X = \prod_{i=0}^{\infty} \{0, 1, \ldots, n_i-1\}$$
  have the product topology. Let
  $U=X-\{\beta_{\max}\}$,
  $V=X-\{{\bf 0}\}$
  and let
  $\varphi\: U \arw V$
  be the ``odometer'' map. Then the crossed product of $C(X)$ by the
partial automorphism
  $\theta\: C_0(U) \arw C_0(V)$,
  induced by $\varphi^{-1}$, is isomorphic to the UHF-algebra
associated to the supernatural number $\prod_{i=0}^{\infty} n_i$.}

  \vskip0pt plus.3\vsize\penalty-250
  \vskip0pt plus -.3\vsize\bigskip\vskip\parskip
  \centerline{\scX References}
  \nobreak\medskip
  \frenchspacing

  \stdbib{\Bratteli}
  {O. Bratteli}
  {Inductive limits of finite dimensional \cstar-algebras}
  {Trans. Amer. Math. Soc.}
  {171}
  {(1972), 195--234}

  \stdbib{\EL}
  {G. A. Elliott and T. Loring}
  {AF embeddings of $C({\bf T}^2)$ with prescribed K-theory}
  {J. Funct. Analysis}
  {103}
  {(1992), 1--25}

  \bib{\NP}
  {R. Exel}
  {Circle actions on \cstar-algebras, partial automorphisms and a
  generalized Pimsner--Voiculescu exact sequence}
  {preprint, University of New Mexico, 1992}

  \stdbib{\Loring}
  {T. Loring}
  {The K-theory of AF embeddings of the rational rotation algebras}
  {K-Theory}
  {4}
  {(1991), 227-243}

  \stdbib{\PV}
  {M. Pimsner and D. Voiculescu}
  {Exact sequences for $K$-groups and Ext-groups of certain
  cross-products \cstar-algebras}
  {J. Oper. Theory}
  {4}
  {(1980), 93-118}

  \stdbib{\PVEmbedd}
  {\sameauthor}
  {Embedding the irrational rotation \cstar-algebra into an
  AF-algebra}
  {J. Oper. Theory}
  {4}
  {(1980), 201--210}

  \bib{\Ped}
  {G. K. Pedersen}
  {\cstar-Algebras and their Automorphism Groups}
  {Academic Press, 1979}

  \bib{\Power}
  {S. C. Power}
  {Limit algebras: An introduction to subalgebras of \cstar-algebras}
  {Pitman Research Notes, to be published by Longman, 1992}

  \stdbib{\Putnam}
  {I. F. Putnam}
  {The \cstar-algebras associated with minimal homeomorphisms of the
  Cantor set}
  {Pacific J. Math.}
  {136}
  {(1989), 329--353}

  \bib{\SV}
  {\c S. \Stratila\ and D. Voiculescu}
  {Representations of AF-algebras and of the group $U(\infty)$}
  {Lecture Notes in Mathematics 486, Springer--Verlag, 1975}

  \bigskip
  \bigskip
  \rightline{November 9, 1992}
  \end